\documentclass[12pt,reqno]{amsart}
\usepackage{amsmath}
\usepackage{amssymb}
\usepackage{latexsym}

\newcommand{\Zint}{\mathbb {Z}}    
     
\newcommand{\Rea}{\mathbb {R}}      % Real number field
\newcommand{\Cplx}{\mathbb {C}}     % Complex  number field

\newcommand{\halmos}{\rule{5pt}{5pt}}

\numberwithin{equation}{section}

\newtheorem{prop}{\bf Proposition}[section]

\newtheorem{lemma}[prop]{\bf Lemma}

\newtheorem{cor}[prop]{\bf Corollary}

\begin{document}
\baselineskip 15pt

\title[Analytic continuation of eigenvalues of the Lam\'e operator]
{Analytic continuation of eigenvalues of the Lam\'e operator}
\author{Kouichi Takemura}
\address{Department of Mathematical Sciences, Yokohama City University, 22-2 Seto, Kanazawa-ku, Yokohama 236-0027, Japan.}
\email{takemura@yokohama-cu.ac.jp}

\begin{abstract}
Eigenvalues of the Lam\'e operator are studied as complex-analytic functions in period $\tau$ of an elliptic function.
We investigate the branching of eigenvalues numerically and clarify the relationship between the branching of eigenvalues and the convergent radius of a perturbation series.
\end{abstract}

\maketitle

\section{Introduction}

The Lam\'e equation is an ordinary differential equation given by
\begin{equation}
 \left( -\frac{d ^2}{d x^2} +n(n+1) \wp (x) \right) f(x) =E f(x),
\label{eqn:lame}
\end{equation}
where $\wp (x)$ is the Weierstrass $\wp$-function which is doubly-periodic with a basic period of $(1,\tau)$, $n\in \Zint _{\geq 1}$ and $E$ is a constant.
In \cite[\S 23]{WW} and \cite[\S 15]{Erd} this equation is discussed in detail.

To analyze the spectral of Eq.(\ref{eqn:lame}), we can choose boundary conditions in various ways.
One is to impose a non-trivial periodic or anti-periodic solution to Eq.(\ref{eqn:lame}). Then, the set of eigenvalues $E$ is discrete and the periodic or the anti-periodic solution is called the Lam\'e function or the singly-periodic Lam\'e function.
Another is to impose a non-trivial doubly-periodic solution to Eq.(\ref{eqn:lame}). In this case the set of eigenvalues $E$ is finite and the doubly-periodic solution is called the Lam\'e polynomial.
When we change the variable by $z=\wp (x)$, the doubly periodic function is essentially expressed as a polynomial in $z$.
Related to quantum mechanics we can choose a boundary condition to have a non-trivial square-integrable solution on the interval $(0,1)$ to Eq.(\ref{eqn:lame}).
We remark that the eigenvalue $E$ with each boundary condition depends on the period $\tau$.

In this paper we investigate how the eigenvalues of Lam\'e functions depend on $\tau$.
In particular we consider branching of the eigenvalues as a complex-analytic function in $\tau $ for the case $n=1$.

Set $q=\exp(\pi \sqrt{-1} \tau)$. 
It is shown in \cite{Tak2} that eigenvalues never stick together if $q \in \Rea$ and $-1< q<1$.
Therefore if $q \in \Rea$ and $-1<q<1$ then there is no branching of the eigenvalue $E$ as a function in $q$ (or $\tau$).

Also note that we can calculate eigenvalues of Lam\'e functions as power series in $q$ by considering perturbation from the trigonometric model (the case $q=0$) as written in \cite{Tak2}.
It is proved in \cite{Tak2} that the convergence radius is not zero.
If the convergence radius is $1$, the eigenvalue is analytic in $\tau $ on the upper half plane, but it is observed numerically that the convergence radii of some eigenvalues are not $1$ (see section \ref{sec:pert}).
Hence there must exist a singularity on the convergence circle.

On the other hand it is known that for the Lam\'e equation with $n\in \Zint _{\geq 1}$ or more generally for the Heun equation with integer coupling constants, the global monodromy is expressed by a hyperelliptic integral \cite{Tak3}. As an application we obtain a condition for $q$ that causes a branching of eigenvalues of the Lam\'e function (see \cite{Tak3} or section \ref{sec:monod} in this paper). 
By thorough calculation, we obtain numerically some values $q$ which produce branching. 

Finally, we find that the absolute value of the branching point calculated by investigating the hyperelliptic integral nearly coincides with the convergence radius calculated by perturbation expansion. In other words we obtain a compatibility between the global monodromy written as a hyperelliptic integral and the perturbation expansion through the branching point.

This paper is organized as follows.
In section \ref{sec:bound} we review several choices for setting the boundary conditions for the Lam\'e operator and observe their relationship.
In section \ref{sec:pert} we explain results on perturbation and the convergence radius.
In section \ref{sec:monod} we consider the global monodromy and search for branching points numerically.
In section \ref{sec:conpa} we discuss the compatibility between perturbation and branching points.
In the appendix, several propositions are proved and definitions and properties of elliptic functions are provided. 

Throughout this paper, we assume that $n$ is a positive integer, and we use the conventions that $f(x)$ is periodic $ \Leftrightarrow f(x+1)=f(x) $, $f(x)$ is anti-periodic $ \Leftrightarrow f(x+1)=-f(x) $ and $f(x)$ is doubly-periodic $ \Leftrightarrow ((f(x+1)=f(x) $ or $-f(x))$ and $(f(x+\tau )=f(x) $ or $-f(x)))$.

\section{Boundary value problems of the Lam\'e operator} \label{sec:bound}

We consider boundary value problems of the Lam\'e operator $H$, where
\begin{align}
& H = -\frac{d ^2}{d x^2} +n(n+1) \wp (x) . \label{LameH}
\end{align}

Let $\sigma _{int} (H)$ be the set of eigenvalues of $H$ whose eigenvector is square-integrable on the interval $(0,1)$, i.e.
\begin{equation}
\sigma _{int} (H) =\{ E | \exists f(x) \in L^2((0,1)) \setminus \{0\}, Hf(x)=Ef(x) \}.
\end{equation}

Let $\sigma _d (H)$ be the set of eigenvalues of $H$ whose eigenvector is doubly-periodic, i.e.
\begin{align}
& \sigma _d (H) = \\
& \{ E | \exists f(x) \neq 0 \mbox{ s.t. }  Hf(x)=Ef(x), \; f(x+1)=\pm f(x) ,  \; f(x+\tau )=\pm f(x) \}, \nonumber 
\end{align}
Note that the doubly-periodic eigenvector is simply the Lam\'e polynomial. It is known \cite{WW} that $\# \sigma _d (H) =2n+1$. 

Let $\sigma _s (H)$ be the set of eigenvalues of $H$ whose eigenvector is singly-periodic.
Set
\begin{align}
& \sigma _p (H) =\{ E | \exists f(x) \neq 0 \mbox{ s.t. }  Hf(x)=Ef(x), f(x+1)=f(x) \}, \\
& \sigma _{ap} (H) =\{ E | \exists f(x) \neq 0 \mbox{ s.t. }  Hf(x)=Ef(x), f(x+1)=-f(x) \}. 
\end{align}
Then $\sigma _s (H) = \sigma _p (H) \coprod \sigma _{ap} (H)$. On the sets $\sigma _{int} (H)$, $\sigma _d (H)$ and $\sigma _s (H)$ we have

\begin{prop} \label{prop:spec1}
(i) For $\tau \in \Rea +\sqrt{-1} \Rea _{>0}$, we have
\begin{eqnarray}
& \sigma _{int} (H) \cup \sigma _d (H) = \sigma _s (H) .&
\end{eqnarray}
(ii) Assume that $q =\exp(\pi \sqrt{-1} \tau) \in \Rea$ and $0<|q|<1$. Then
\begin{eqnarray}
& \sigma _{int} (H) \coprod \sigma _d (H) = \sigma _s (H) ,& \label{eq:disj}
\end{eqnarray}
i.e., $\sigma _{int} (H) \cup \sigma _d (H) = \sigma _s (H)$ and $\sigma _{int} (H) \cap \sigma _d (H) = \phi$.
\end{prop}
We prove this proposition in the appendix.
Note that, if $q$ is not real, then the proposition $\sigma _{int} (H) \cap \sigma _d (H) = \phi$ might be false. In fact, if $n=1$ and $q=\sqrt{-1}(.3281\dots )$, then it seems that $-e_1 \in \sigma _{int} (H) \cap \sigma _d (H)$ (see Proposition \ref{prop:01} and Table 3). 

Next, we briefly explain the relationship to the finite-gap potential.
Let 
\begin{align}
& I = -\frac{d ^2}{d x^2} +n(n+1) \wp (x+\tau /2) 
\end{align}
and $\sigma _b(I)$ be the set such that 
$$
E \in \sigma _b(I) \; \Leftrightarrow \mbox{ Every solution to }(I-E)f(x)=0 \mbox{ is bounded on }x \in \Rea .
$$
Ince \cite{I} established that, if $q =\exp(\pi \sqrt{-1} \tau) \in \Rea$, then  
\begin{equation}
\Rea \setminus \overline{\sigma _b(H)}= (-\infty,E_{0})\cup (E_{1},E_{2})\cup \dots \cup (E_{2n-1}, E_{2n})
\end{equation}
where $\overline{\sigma _b(H)}$ is the closure of the set $\sigma _b(H)$ in $\Cplx $, $E_i \in \sigma _d (H)$ and $E_0<E_{1}<\cdots <E_{2n}$.
Hence there is a finite band structure on eigenvalues of unbounded eigenvectors.
This is referred to as finite-band potential or finite-gap potential.

\section{Perturbation and convergence radius} \label{sec:pert}

In this section we calculate eigenvalues of Lam\'e functions as power series in $q (= \exp (\pi \sqrt{-1} \tau ))$.
For this purpose we consider perturbation from the trigonometric model. First we consider a trigonometric limit $q \rightarrow 0 \: (\Leftrightarrow \tau \rightarrow \sqrt{-1} \infty)$ and later apply a method of perturbation from the trigonometric model.

For the case $q=0$ the spectral problem becomes much simpler. Set 
\begin{align}
& H_T=  -\frac{d^2}{dx^2} + n(n+1)\frac{\pi ^2}{\sin ^2\pi x} .\label{triIno} 
\end{align}
Then $H \rightarrow H_T -\frac{\pi^2}{3} n (n +1)$ as $q = \exp ( \pi \sqrt{-1} \tau )\rightarrow 0$.
The operator $H_T$ is the Hamiltonian of the P\"oschl-Teller system or the $A_1$ trigonometric Calogero-Moser-Sutherland system.
Set
\begin{equation}
\Phi(x)=(\sin \pi x)^{n+1}, \; \; \; v_m=  \tilde{c}_m C^{n+1}_m (\cos \pi x) \Phi(x) ,\; \; \; (m \in \Zint_{\geq 0}),
\label{Gegenpol}
\end{equation}
where the function $C^{\nu}_m (z)=\frac{\Gamma (m+2\nu)}{m! \Gamma(2\nu)}\, _2 \! F_1(-m,m+2\nu ; \nu+\frac{1}{2} ;\frac{1-z}{2})$ is the Gegenbauer polynomial of degree $m$ and $\tilde{c} _m=\sqrt{\frac{2^{2n+1}(m+n+1)m! \Gamma (n+1)^2}{\Gamma (m+2n+2)}}$.
Then
\begin{align}
H_T v_m=\pi^2(m+n+1)^2 v_m,
\end{align}
and $\langle v_m, v_{m'} \rangle =\delta_{m,m'}$, where the inner product is defined  by 
\begin{equation}
\langle f,g\rangle =\int_{0}^{1} \overline{f(x)} g(x) dx.
% \; \; \langle f,g\rangle _{\Phi} =\int_{0}^{1} dx\overline{f(x)} g(x) |\Phi(x)|^2.
\label{innerprod}
\end{equation}

Set
\begin{equation}
{\bold H} = \left\{ f \! : \Rea \rightarrow \Cplx \mbox{ measurable} \left|
\begin{array}{l}
 \int_{0}^{1} |f(x)| ^2 dx<+\infty, \\
 f(x)= f(x+2) \mbox{ a.e. }x, \\
 f(x)=(-1)^{n+1} f(-x) \mbox{ a.e. }x
\end{array}
\right. \right\} ,
\label{Hilb1l10}
\end{equation}
\begin{align}
& {\bold H}_+ = \{ f \in {\bold H} | f(x)=f(x+1) \mbox{ a.e. }x  \} \nonumber ,\\
& {\bold H}_- = \{ f \in {\bold H} | f(x)=-f(x+1) \mbox{ a.e. }x  \} \nonumber .
\end{align}
Inner products on the Hilbert space ${\bold H}$ and its subspaces ${\bold H}_+$, ${\bold H}_-$ are given by $\langle \cdot , \cdot  \rangle $. Then we have ${\bold H}_+ \perp {\bold H}_-$ and ${\bold H}= {\bold H}_+ \oplus {\bold H}_-$.
The Hamiltonian $H$ (see Eq.(\ref{LameH})) acts on a certain dense subspace of ${\bold H}$ (resp. ${\bold H}_+$, ${\bold H}_-$) and the space spanned by functions $\{ v _m | m\in \Zint_{\geq 0} \}$ (resp. $\{ v _m | m\in 2 \Zint_{\geq 0} \}$, $\{ v _m | m\in 2\Zint_{\geq 0} +1 \}$) is dense in ${\bold H}$ (resp. ${\bold H}_+$, ${\bold H}_-$).

Now we apply a method of perturbation and have an algorithm for obtaining eigenvalues and eigenfunctions as formal power series of $q$. For details see \cite{Tak2}.

Set $q=\exp(\pi \sqrt{-1} \tau) $. For the Lam\'e operator (see Eq.(\ref{LameH})), we adopt the notation $H(q)$ instead of $H$. The operator $H(q)$ admits the following expansion:
\begin{equation}
H(q)(=H)=H_T  -\frac{\pi^2}{3} n (n +1)+\sum_{k=1}^{\infty} V_{2k}(x) q^{2k},
\label{Hamilp0}
\end{equation}
where $H_T$ is the Hamiltonian of the trigonometric model and $V_{2k}(x)$ are functions in $x$ which are determined by using Eq.(\ref{wpth}).

Set 
\begin{equation}
E_m=\pi^2(m+n+1)^2 -\frac{\pi^2}{3} n (n +1) .
\end{equation}
Then $v_m$ is an eigenfunction of the operator $H(0)$ with the eigenvalue $E_m$.

Based on the eigenvalues $E_m$ $(m \in \Zint _{\geq 0} )$ and the eigenfunctions $v_m$ of the operator $H(0)$, we determine eigenvalues $E_m (q) = E_m+\sum_{k=1}^{\infty} E_m^{\{2k\}}q^{2k}$ and normalized eigenfunctions $v_m(q)= v_m+ \sum_{k=1}^{\infty} \sum_{m'\in \Zint _{\geq 0}} c_{m ,m'}^{\{ 2k \} }v_{m'}q^{2k}$ of the operator $H(q)$
as formal power series in $q$.
In other words, we will find $E_m(q) $ and $v_m(q)$ that satisfy equations
\begin{eqnarray}
& H(q)v_m(q)=(H(0) +\sum_{k=1}^{\infty} V_{2k}(x) q^{2k} )v_m(q) = E_m(q)v_m(q) ,\label{Hpertexp} &\\
& \langle v_m(q) ,  v_m(q) \rangle =1 , & \nonumber
\end{eqnarray}
as formal power series of $q$.

First we calculate coefficients $\sum_{m' \in \Zint _{\geq 0}} d_{m ,m'}^{\{2k\}}v_{m'}= V_{2k}(x) v_m$ $(k\in \Zint_{>0}, \; m \in \Zint _{\geq 0})$. 
Next we compute $E_m^{\{2k\}}$ and $ c_{m ,m'}^{\{2k\}}$ for $k \geq 1$ and $m , m'\in \Zint _{\geq 0}$.  By comparing coefficients of $v_{m'}q^{2k}$ in Eq.(\ref{Hpertexp}), we obtain recursive relations for $E_m^{\{2k\}}$ and $c_{m ,m'}^{\{2k\}}$. 
For details see \cite{Tak2}.
Note that, if $m-m'$ is odd, then $d_{m ,m'}^{\{2k\}}=c_{m ,m'}^{\{2k\}}=0$.
Convergence of the formal power series of eigenvalues in the variable $q$ obtained by the algorithm of perturbation is shown in \cite{Tak2}.
\begin{prop} \label{cor:conv} \cite[Corollary 3.7]{Tak2}
Let $E_m(q)$ $(m\in \Zint _{\geq 0})$ (resp. $v_m(q)$) be the formal eigenvalue (resp. eigenfunction) of the Hamiltonian $H(q)$ defined by Eq.(\ref{Hpertexp}).
If $|q|$ is sufficiently small then the power series $E_m(q)$ converges and as an element in the Hilbert space ${\bold H}$ the power series $v_m(q)$ converges.
\end{prop}

We show an expansion of the first few terms of the eigenvalue $E_m(q)$ and the radius of convergence for the case $n=1$ in Table 1.
We calculate the expansion of $E_m(q) =E_m +\sum _{k} E_m ^{\{2k\} }q^{2k}$ for more than $100$ terms and approximate the absolute values of coefficients $ E_m ^{\{2k\} }$ by $a b^{2k}$ for some constants $a$ and $b$ which are determined by the method of least squares. Then, the radius of convergence is inferred by $\liminf _{k \rightarrow \infty} 1/( |E_m ^{\{2k\} }|/a)^{1/2k}$. 
The inferred radius of convergence and expansions of the first few terms of the eigenvalue $E_m(q)$ are calculated as follows:
%\begin{scriptsize}
\begin{center}
\begin{tabular}{|l|l|l|} \hline
$E_0 (q)$ & $\pi^2 \left( \frac{10}{3} + \frac{80}{3} q^2+\frac{1360}{27} q^4+\frac{20800}{243} q^6 +\frac{195920}{2187}q^8+ \frac{3174880}{19683} q^{10} +\frac{684960}{59049}q^{12}+ \dots \right)$ & .749 \\ \hline
$E_2 (q)$ & $\pi^2 \left( \frac{46}{3} + \frac{272}{15} q^2+\frac{198928}{3375} q^4+\frac{55403584}{759375} q^6 +\frac{4307155408}{34171875}q^8 + \frac{2879355070048}{38443359375} q^{10} +\dots \right)$ & .749 \\ \hline
$E_4 (q)$ & $\pi^2 \left( \frac{106}{3} + \frac{592}{35} q^2 +\frac{2279248}{42875} q^4+\frac{3773733184}{52521875} q^6 +\frac{1634762851088}{12867859375}q^8 +\dots \right)$ & .875 \\ \hline \hline
$E_1 (q)$ & $\pi^2 \left( \frac{25}{3} + 20 q^2+ 65 q^4+\frac{115}{2} q^6 +\frac{2165}{16}q^8+ \frac{3165}{32} q^{10} +\frac{23965}{128}q^{12} +\frac{38755}{256}q^{14} +\dots \right)$ & .838 \\ \hline
$E_3 (q)$ & $\pi^2 \left( \frac{73}{3} + \frac{52}{3} q^2 +\frac{1493}{27} q^4+\frac{35671}{486} q^6 +\frac{4492153}{34992}q^8 + \frac{55853449}{629856} q^{10} +\frac{1646085467}{7558272} q^{12} + \dots \right)$ & .838 \\ \hline
$E_5 (q)$ & $\pi^2 \left( \frac{241}{3} + \frac{82}{5} q^2+\frac{50339}{1000} q^4+\frac{13640101}{200000} q^6 +\frac{3872868499}{32000000}q^8+ \frac{3267409458867}{32000000000} q^{10} +\dots \right)$ & .906 \\ \hline
\end{tabular} 
\end{center}
%\end{scriptsize}
\noindent {\it Table 1. Expansion of the first few terms and the inferred radius of convergence.}

We introduce propositions on the spectral of the Hamiltonian $H$ on the Hilbert spaces for the case $q^2 \in \Rea$ and $|q|<1$.
Let $\sigma _{{\bf H}}(H)$ (resp. $\sigma _{{\bf H}_+}(H)$, $\sigma _{{\bf H}_-}(H)$) be the spectral of the operator $H$ on the space ${\bf H}$ (resp. ${\bf H}_+$, ${\bf H}_-$).

\begin{prop} \label{prop:discr} (c.f. \cite[Propositions 3.2, 3.5]{Tak2})
Let $q^2 \in \Rea$ and $|q|<1$. The operator $H$ is essentially selfadjoint on the Hilbert space ${\bf H}$ (resp. ${\bf H}_+$, ${\bf H}_-$). The spectrum $\sigma _{{\bf H}}(H)$ (resp. $\sigma _{{\bf H}_+}(H)$, $\sigma _{{\bf H}_-}(H)$) contains only point spectra and it is discrete.
\end{prop}

\begin{prop} \label{mainthmKato} (c.f. \cite[Theorem 3.6]{Tak2})
Let $q^2 \in \Rea$ and $|q|<1$. All eigenvalues of $H$ on the space ${\bf H}$ can be represented as $E_m (q)$ $(m \in \Zint _{\geq 0})$, which is real-holomorphic in $q^2 \in (-1,1)$ and $E_m (0) =E_m$.
The eigenfunction $v_m (q)$ of the eigenvalue $E_m (q)$ is holomorphic in $q^2 \in (-1,1)$ as an element in $L^2$-space, and the eigenvectors $v_m (q)$ $(m \in \Zint _{\geq 0})$ form a complete orthonormal family on ${\bf H}$.
\end{prop}

It is shown that, if $q^2 \in \Rea$, $|q|<1$ and $m \in 2\Zint _{\geq 0}$ (resp. $m \in 2\Zint _{\geq 0} +1$), then the corresponding eigenvector $v_m (q)$ belongs to the space ${\bf H}_+$ (resp. ${\bf H}_-$) and we have
\begin{align}
& \sigma _{{\bf H}}(H)= \{ E_m (q)| m \in \Zint _{\geq 0} \} \\
& \sigma _{{\bf H}_+}(H)= \{ E_m (q)| m \in 2\Zint _{\geq 0} \} \nonumber \\
& \sigma _{{\bf H}_-}(H)= \{ E_m (q)| m \in 2\Zint _{\geq 0} +1 \} \nonumber 
\end{align}

Among the spaces $\sigma _{{\bf H}} (H)$, $\sigma _{{\bf H}_+}(H)$, $\sigma _{{\bf H}_-}(H)$, $\sigma _{int} (H)$, $\sigma _{p}(H)$ and $\sigma _{ap}(H)$, the following relations are satisfied:
\begin{prop} \label{prop:spec2}
We have $\sigma _{{\bf H}} (H) = \sigma _{int} (H)$, $\sigma _{{\bf H}_+}(H) =\sigma _{int} (H) \cap \sigma _p(H)$ and $\sigma _{{\bf H}_-}(H) =\sigma _{int} (H) \cap \sigma _{ap}(H)$
\end{prop}
\begin{proof}
It follows from the definition of ${\bold H}$ that, if $f(x) \in {\bold H}$, then the function $f (x)$ is square-integrable on $(0,1)$, i.e. $\sigma _{{\bold H}} (H) \subset \sigma _{int} (H)$. Now we show $\sigma _{int} (H) \subset \sigma _{{\bold H}} (H)$. 
Let $E \in \sigma _{int} (H)$. Then there exists a non-zero function $f(x)$ such that $Hf(x)=Ef(x)$ and $\int _0^1 |f(x)|^2 dx < \infty$. The exponent of the differential equation $(H-E)f(x)=0$ at $x=0$ is $\{ -n, n+1\} $.
Since the function $f(x)$ is square-integrable and the equation $(H-E)f(x)=0$ is invariant under the transformation $x \leftrightarrow -x$, the function $f(x)$ is expanded as
\begin{equation}
f(x)= x^{n+1} (c_0 + c_1 x^2 +c_2 x^4 +\dots ) \; \; \; (c_0 \neq 0)
\end{equation}
and satisfies $f(x) =(-1)^{n+1} f(-x)$. From the periodicity, the function $f(x+1)$ is also an eigenfunction. The function $f(x+1)$ is written as a linear combination of $f(x)$ and another linearly independent solution, and we have $f(x+1)= C f(x)$ for some $C(\neq 0)$ because $f(x)$ is locally square-integrable near $x=1$. It follows immediately that $f(-x-1)=C^{-1} f(-x)$. We have $C f(x)= f(x+1)=(-1)^{n+1} f(-x-1) = (-1)^{n+1} C^{-1} f(-x) = C^{-1} f(x)$. Hence $C \in \{ \pm 1\}$ and $f(x+2)=f(x)$.
Therefore we have $f(x) \in {\bold H}$, $E \in \sigma _{{\bold H}} (H)$ and $\sigma _{int} (H) \subset \sigma _{{\bold H}} (H)$.

Relations $\sigma _{{\bf H}_+}(H) =\sigma _{int} (H) \cap \sigma _p(H)$ and $\sigma _{{\bf H}_-}(H) =\sigma _{int} (H) \cap \sigma _{ap}(H)$ are obtained by considering periodicity.
\end{proof}

It is shown that eigenvalues never stick together as in \cite{Tak2}.

\begin{prop} \label{thm:multone} (c.f. \cite[Theorem 3.9]{Tak2})
Let $E_m(q)$ $(m\in \Zint_{\geq 0})$ be the eigenvalues of $H(q)$ defined in Proposition \ref{mainthmKato}. If $q^2 \in \Rea$ and $|q|<1$, then $E_m(q)\neq E_{m'}(q)$ $(m\neq m')$. In other words, eigenvalues never stick together under the condition $q^2 \in \Rea$ and $|q|<1$.
\end{prop}
\begin{proof}
Assume that the proposition is wrong. Then there exists $m$ and $q$ such that $E_m(q)=E_{m+1}(q)$. Let $f(x)$ and $\tilde{f}(x)$ be the corresponding eigenfunctions. Then one of $f(x)$ or $\tilde{f}(x)$ is periodic and the other is anti-periodic. Hence $f(x)$ and $\tilde{f}(x)$ are linearly independent. Since there is no first differential term in $H$, we have $\left( \frac{d^2}{dx^2}f(x) \right) \tilde{f}(x)- f(x) \frac{d^2}{dx^2}\tilde{f}(x) =0$. Hence $\left( \frac{d}{dx}f(x) \right) \tilde{f}(x)- f(x) \frac{d}{dx}\tilde{f}(x)$ is a constant and it is non-zero by linear independence. It contradicts the periodicity of $f(x)$ and $\tilde{f}(x)$ and we obtain the proposition.
\end{proof}

\begin{cor} \label{cor:ineq} (c.f. \cite[Corollary 3.10]{Tak2})
If $q^2 \in \Rea$, $|q|<1$ and $m<m'$, then $E_m(q)< E_{m'}(q)$.
\end{cor}

\section{Monodromy and branching points} \label{sec:monod}

We consider the monodromy of solutions of 
\begin{equation}
 H f(x) =E f(x),\; \; H= -\frac{d ^2}{d x^2} +2 \wp (x)
\label{eqn:lame1}
\end{equation}
for each $E$. Note that this is the case $n=1$ in Eq.(\ref{eqn:lame}). 

For the case $n=1$, we have $\sigma _d (H)= \{ -e_1, -e_2, -e_3 \}$ and the corresponding doubly-periodic eigenfunctions are $\wp_1(x), \wp_2(x), \wp_3(x)$ (see Eq.(\ref{wpi})). From the periodicity of $\wp_i(x)$ $(i=1,2,3)$ we have $\sigma _d (H) \cap \sigma _p (H)= \{ -e_1 \}$ and $\sigma _d (H)  \cap \sigma _{ap} (H)= \{ -e_2, -e_3 \}$.

We now consider the expression of solutions to Eq.(\ref{eqn:lame1}) for each $E$.
The functions $\Xi(x,E)$ and $P(E)$ defined around Proposition \ref{prop:prod} for the case $n=1$ are calculated as $\Xi(x,E) = \wp \left(x \right)+E$ and $P(E)=(E+e_1)(E+e_2)(E+e_3)$. Then the function $\Lambda ( x, E)$ defined in Eq.(\ref{integ1}) is a solution to the differential equation (\ref{eqn:lame}) (see Proposition \ref{prop:Lambda}), and
it is also expressed as 
\begin{equation}
\Lambda (x, E) = A  \frac{\sigma (x+t_0)}{\sigma (x)}e^{-x\zeta (t_0)} , \quad E=-\wp (t_0),
\label{eq:Lamt0}
\end{equation}
for suitably chosen $A$ (see \cite[\S 39]{P} or \cite[\S 23.7]{WW}), where $\sigma (x)$ is the Weierstrass sigma-function and $\zeta (x)$ is the Weierstrass zeta-function (see Appendix). Note that we can show directly that the function $\Lambda (x, E)$ written as Eq.(\ref{eq:Lamt0}) satisfies Eq.(\ref{eqn:lame1}).
It follows from Eq.(\ref{eq:Lamt0}) and Eq.(\ref{eq:sigmaper}) that the monodromy is described as 
\begin{equation} 
\Lambda (x+1,E)=\Lambda (x,E) \exp \left( 2\eta _1 t_0-\zeta (t_0) \right),
\label{analcont0001}
\end{equation}
where $\eta _1=\zeta (1/2)$.
Hence, if $ 2\eta _1 t_0-\zeta (t_0) \in \pi \sqrt{-1} \Zint $ (resp. $2\eta _1 t_0-\zeta (t_0) \in 2\pi \sqrt{-1} \Zint $, $2\eta _1 t_0-\zeta (t_0) \in 2\pi \sqrt{-1} \Zint + \pi \sqrt{-1}$), then $E \in \sigma _s (H)$ (resp. $E \in \sigma _p (H)$, $E \in \sigma _{ap} (H)$).
It follows from Proposition \ref{prop:Tak344} that, if $2\eta _1 t_0-\zeta (t_0) \not\in \pi \sqrt{-1} \Zint $, then $E \not\in \sigma _{int} (H)$.
By Proposition \ref{prop:spec1} and Proposition \ref{prop:spec2}, if $-1<q(=\exp(\pi \sqrt{-1} \tau ))<1$, then we have 
\begin{align}
& \sigma _{{\bf H}}(H) = \sigma _{s}(H) \setminus \{ -e_1, -e_2, -e_3 \} ,\\
& \sigma _{{\bf H}_+}(H) = \sigma _{p}(H) \setminus \{ -e_1 \} , \; \; \; \sigma _{{\bf H}_-}(H) = \sigma _{ap}(H) \setminus \{ -e_2, -e_3 \} . \nonumber
\end{align}

The eigenvalue in $\sigma _p (H)$ is analytically continued in $q$ (or $\tau$) as to preserve the property 
\begin{equation}
 E = -\wp(t_0), \; \; \;  2\eta _1 t_0-\zeta (t_0)  \in 2\pi \sqrt{-1} \Zint.
\label{analcontcond}
\end{equation}
and the eigenvalue in $\sigma _{ap} (H)$ is analytically continued in $q$ (or $\tau$) as to preserve the property 
\begin{equation}
 E = -\wp(t_0), \; \; \;   2\eta _1 t_0-\zeta (t_0) \in 2\pi \sqrt{-1} \Zint + \pi \sqrt{-1}.
\label{analcontcondap}
\end{equation}
It follows from the relation $E=-\wp (t_0)$ and Eq.(\ref{ellinteg1}) that Eq.(\ref{analcont0001}) is rewritten as 
\begin{equation} 
\Lambda (x+1,E)=\Lambda (x,E) \exp \left( -\frac{1}{2}\int_{-e_1}^{E}\frac{\tilde{E}-2\eta_1}{\sqrt{-(\tilde{E}+e_1)(\tilde{E}+e_2)(\tilde{E}+e_3)}} d\tilde{E}\right),
\label{analcont0001int}
\end{equation}
Hence we reproduce the monodromy formula in terms of (hyper)elliptic integral which was obtained in \cite{Tak3}.
For analyticity of elements in $\sigma _{p} (H)$ or $\sigma _{ap} (H)$, we have
\begin{prop} (c.f. \cite[Theorem 4.6 (ii)]{Tak3}) \label{prop:00}
If the eigenvalue $E$ satisfies Eq.(\ref{analcontcond}) or Eq.(\ref{analcontcondap}), $E-2\eta_1\neq 0$ and $E \neq -e_1, -e_2, -e_3$ at $q=q_*$, then the eigenvalue $E$ satisfying Eq.(\ref{analcontcond}) or Eq.(\ref{analcontcondap}) is analytic in $q$ around $q=q_*$.
\end{prop}
Note that Proposition \ref{prop:00} is proved by applying the implicit function theorem as is done in \cite[Theorem 4.6 (ii)]{Tak3}. The following proposition describes the condition for $q$ (or $\tau $) that the set $\sigma _{d} (H) \cap \sigma _{int} (H)$ is non-empty.
\begin{prop} \label{prop:01}
Under the assumption $E \in \sigma _{d} (H)$ (i.e., $E \in \{ -e_1, -e_2, -e_3\} $), the condition $E \in \sigma _{int} (H)$ is equivalent to the condition $E-2\eta _1=0$.
\end{prop}
\begin{proof}
It follows from the assumption that $E=-e_i$ for some $i\in \{1,2,3\}$.
A solution to Eq.(\ref{eqn:lame1}) for $E=-e_i$ is written as $\wp _i (x)$, and another solution is written as $\wp _i (x) \int  (1/\wp _i (x) ^2) dx$. By Eqs.(\ref{eq:sigmaper}, \ref{sigmai}) we have
\begin{equation}
\int  \frac{dx}{\wp _i (x) ^2} = \int  \frac{dx}{\wp (x)-e_i} = \int  \frac{(\wp (x+\omega _i)-e_i)dx}{(e_i-e_{i'})(e_i-e_{i''})} = -\frac{\zeta (x+\omega _i)+e_i x}{(e_i-e_{i'})(e_i-e_{i''})},
\end{equation}
where $i', i'' \in \{1,2,3\}$ with $i'<i''$, $i\neq i'$, and $i\neq i''$.
Set $s_1 (x)= \wp _i (x)$ and $s_2 (x)= \wp _i (x) (\zeta (x+\omega _i)+e_i x -\eta _i) $. Then they are a basis of solutions to Eq.(\ref{eqn:lame1}) for $E=-e_i$, and $s_1 (x)$ (resp. $s_2 (x)$) is odd (resp. even).
Since $s_1 (x)$ has a pole at $x=0$ and $s_2 (x)$ is holomorphic at $x=0$, square-integrable eigenfunction on $(0,1)$ is written as $A s_2 (x)$ for some constant $A$.
Since $s_2 (x+1)$ cannot have a pole at $x=0$ for square-integrability and it is written as 
\begin{align}
s_2 (x+1)& =  \wp _i(x+1) \left( \zeta (x+\omega _i+1)+(x+1)e_i -\eta _i \right) \label{eq:s2}\\
& = \pm \left( s_2 (x) +(e_i+2\eta_1) \wp _i (x) \right)
\nonumber
\end{align}
for some sign $\pm$, we have $E-2\eta_1=0$ (i.e., $-e_i-2\eta_1=0$).

Conversely, if $E-2\eta_1=0$ and $E=-e_i$, then it follows from Eq.(\ref{eq:s2}) that $s_2(x)$ is perioic with a period $1$ and it is holomorphic on $\Rea $. Hence $s_2(x)$ is square-integrable on $(0,1)$, and we have $E \in \sigma _{int} (H)$.
\end{proof}

By Propositions \ref{prop:spec1}, \ref{prop:00} and \ref{prop:01}, it follows that if the eigenvalue $E$ in $\sigma _{p} (H)$ or $\sigma _{ap} (H)$ has a branching at $q$, then we have $E-2\eta_1 = 0$.
%The condition $E-2\eta_1 = 0$ is a necessary condition that the eigenvalue $E$ satisfying Eq.(\ref{analcontcond}) or Eq.(\ref{analcontcondap}) has a branching.
Hence a necessary condition that the eigenvalue $E$ in $\sigma _{p} (H)$ or $\sigma _{ap} (H)$ has a branching is that $q$ and $t_0$ satisfy the following conditions:
\begin{align}
& 2\eta_1 = -\wp(t_0) (=E), \label{acc21} \\
& 2 \eta _1 t_0-\zeta (t_0) \in \pi \sqrt{-1} \Zint. \label{acc22}
\end{align}
We try to solve Eqs.(\ref{acc21}, \ref{acc22}) numerically. 
First we fix the value $q$. We expand $\eta_1$, $\wp(t_0)$ and $\zeta (t_0)$ in $q$ according to Eq.(\ref{wpth}) with approximately 100 terms, and solve Eq.(\ref{acc21}) numerically by Newton's method and obtain $t_0$.
We evaluate Eq.(\ref{acc22}) using $t_0$ and check whether it is satisfied or not. Note that the imaginary part of the value $t_0$ should be taken to be small in order to exhibit good convergence.

By investigating more than 1000 complex numbers which satisfy $|q|<.90$, $\Re q \geq 0$ and $\Im q \geq 0$ where $\Re q $ (resp. $\Im q$) is the real part (resp. the imaginary part) of the number $q$, we obtain numerically that the numbers in Table 2 may have branches (i.e. they satisfy Eq.(\ref{acc21}) and Eq.(\ref{acc22})). Note that it seems some numbers do not generate branching.
%\begin{small}
\begin{center}
\begin{tabular}{|l|ll|}
\hline
periodic & $q=$ & $.328106 I, .258666+.697448 I, .510303+.546057 I$ \\
     & &$.746852+.452463 I, .224582 +.842777 I, .552288 +.677536 I$ \\
     & &$.314813+.821858 I, .686317+.559106 I $ \\
\hline
anti-periodic & $q=$ & $.281417+.534362 I, .655163+.503275 I, .264829+.792687 I $ \\
     & & $.535905+.640487 I, .807197+.405705 I $\\
\hline
\end{tabular}
\end{center}
%\end{small}
\noindent {\it Table 2. Numbers which may have branches.}

Next we consider how to continue the eigenvalues analytically in $q$ along a path. Let ${\mathcal C}$ be a path in the complex plane.
The eigenvalue $E$ is continued analytically in $q$ along the path ${\mathcal C}$ by keeping the conditions 
\begin{align}
& E = -\wp(t_0), \label{acc23} \\
& \exists m \in \Zint, \; \; \; 2 \eta _1 t_0 -\zeta (t_0) = m \pi \sqrt{-1}. \label{acc24}
\end{align}
Note that the eigenvalue satisfying Eq.(\ref{acc23}) and  Eq.(\ref{acc24}) for $m\in 2\Zint $ (resp. $m\in 2\Zint +1$) is continued from the eigenvalue in ${\bf H}_+$ (resp. ${\bf H}_-$).

We solve Eqs.(\ref{acc23}, \ref{acc24}) for points which are selected appropriately on the path ${\mathcal C}$ and are connected by choosing close solutions. Note that for each $E$ and $q$ satisfying Eqs.(\ref{acc23}, \ref{acc24}), solutions $(t_0, m)$ may not be unique. Sometimes we need to change to another solution $(t'_0, m')$ to avoid the divergence of continued solutions in $q$.

We continue the eigenvalue $E$ analytically around the possible branches in Table 2. We obtain that the following numbers would not cause branching and they all would satisify $2\eta _1 =- e_i$ for some $i \in \{ 1,2,3 \}$:
%\begin{small}
\begin{center}
\begin{tabular}{|l|l||l|l|}
\hline
$ q= .328106 I $& $2\eta _1 =- e_1$ &  $q= .281417+.534362 I$ & $2\eta _1 =- e_2$ \\
$ q= .510303+.546057 I $ & $2\eta _1 =- e_1$ & $q= .655163+.503275 I$ & $2\eta _1 =- e_2$ \\
$ q= .746852+.452463 I$ & $2\eta _1 =- e_1$ &  $q= .264829+.792687 I$ & $2\eta _1 =- e_3$ \\
  & & $q= .807197+.405705 I$ & $2\eta _1 =- e_2$ \\
\hline
\end{tabular}
\end{center}
%\end{small}
\noindent {\it Table 3. Numbers that do not cause branching.}

For these cases, it is inferred from Proposition \ref{prop:01} that one of the eigenvalues $E_m(q)$ ($m\in \Zint_{\geq 0})$ meets with an eigevalue with doubly-periodic eigenfunction (i.e. $-e_1$, $-e_2$ or $-e_3$).

Let $a \in \Cplx $ and ${\mathcal C}_a$ be the cycle starting from $\Re a$, approaching the point $a$ parallel to the imaginary axis, turning anti-clockwise around $a$ and returning to $\Re a$ as shown in Figure 4.
\begin{center}
\begin{picture}(200,100)(0,0)
\put(10,10){\vector(1,0){170}}
\put(10,10){\vector(0,1){85}}
\put(112,10){\vector(0,1){40}}
\put(112,50){\line (0,1){19}}
\put(108,69){\vector(0,-1){40}}
\put(108,30){\line(0,-1){20}}
\put(110,75){\circle*{2}}
\put(108,75){\oval(10,12)[l]}
\put(112,75){\oval(10,12)[r]}
\put(112,81){\vector(-1,0){4}}
\put(185,6){Re}
\put(15,85){Im}
\put(108,1){$\Re a$}
\put(118,68){$a$}
\put(96,45){${\mathcal C}_a$}
\end{picture}
\end{center}
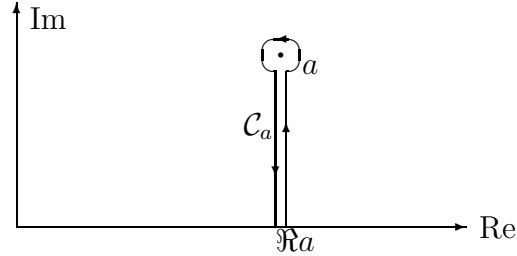
\noindent {\it Figure 4. Cycle ${\mathcal C}_a$.}

We continue the eigenvalue $E$ analytically along the cycle ${\mathcal C}_a$ where $a$ is a branching point which is listed in Table 2 and not listed in Table 3. The branching along the cycle ${\mathcal C}_a$ is then determined as shown in Table 5.
%\begin{small}
\begin{center}
\begin{tabular}{|l|l|}
\hline
$a=.258666+.697448 I$ & $E_0 (q)  \Rightarrow E_2 (q), \; E_2 (q) \Rightarrow E_0 (q), \: \; E_4 (q) \Rightarrow E_4 (q), \; E_6 (q) \Rightarrow E_6 (q)$ \\
$a=.224582 +.842777 I$ & $E_0 (q) \Rightarrow E_4 (q), \; E_2 (q) \Rightarrow E_2 (q), \: \; E_4 (q) \Rightarrow E_0 (q), \; E_6 (q) \Rightarrow E_6 (q)$ \\
$a=.552288 +.677536 I$ & $E_0 (q) \Rightarrow E_4 (q), \; E_2 (q) \Rightarrow E_2 (q), \: \; E_4 (q) \Rightarrow E_0 (q), \; E_6 (q) \Rightarrow E_6 (q)$ \\
$a=.314813+.821858 I$ & $E_0 (q) \Rightarrow E_4 (q), \; E_2 (q) \Rightarrow E_2 (q), \: \; E_4 (q) \Rightarrow E_0 (q), \; E_6 (q) \Rightarrow E_6 (q)$ \\
$a=.686317+.559106 I$ & $E_0 (q) \Rightarrow E_0 (q), \; E_2 (q) \Rightarrow E_4 (q), \: \; E_4 (q) \Rightarrow E_2 (q), \; E_6 (q) \Rightarrow E_6 (q)$ \\
\hline
$a=.535905+.640487 I$ & $E_1 (q) \Rightarrow E_3(q) , \; E_3(q)  \Rightarrow E_1(q) , \: \; E_5(q)  \Rightarrow E_5(q) , \; E_7(q)  \Rightarrow E_7(q) $ \\
\hline
\end{tabular}
\end{center}
%\end{small}
\noindent {\it Table 5. Branching along the cycle ${\mathcal C}_a$}

\section{Convergence radius and branching points} \label{sec:conpa}

In section \ref{sec:monod} we calculated the positions of the branching points of the eigenvalues $E_m(q)$ ($m\in \Zint $) in $q$ and described how the eigenvalues are continued along cycles. In this section we observe that the convergence radii of the eigenvalues $E_m(q)$ calculated by perturbation are compatible with the positions of the branching points.

For the periodic case the closest branching point from the origin is $q=.258666+.697448 I$ ($|q|=.743869$) and the eigenvalues $E_0(q)$ and $E_2(q)$ are connected by continuing analytically along the cycle ${\mathcal C}_q$ $(q=.258666+.697448I)$ (see Table 5). 
It is known that the convergence radius of a complex function expanded at an origin is equal to the distance from the origin to the closest singular point. Hence the convergence radii of the eigenvalues $E_0(q)$ and $E_2(q)$ are both $.743869$.

On the other hand in section \ref{sec:pert} we obtained that the convergence radii of the expansions of the eigenvalues $E_0(q)$ and $E_2(q)$ around $q=0$, calculated by the method of perturbation, are both around $.749$.

Thus, convergence radii calculated by different methods are very close and compatibility between the method of perturbation and the method of monodromy is confirmed. Moreover, we obtain a reason why the convergence radii of the eigenvalues $E_0(q)$ and $E_2(q)$ calculated in section \ref{sec:pert} are very close by considering the branching point.
To get more precise values of convergence radii calculated by perturbation, it is necessary to calculate more terms in $k$ on the expansion $E_m(q)= E_m +\sum _{k} E_m ^{\{2k\} }q^{2k}$ $(m=0,2)$. Generally speaking, it would be impractical to guess a convergence radius numerically from Taylor's expansion.

The second closest branching point from the origin for the periodic case is $q=.224582 +.842777 I$ ($|q|=.872187$) and the eigenvalues $E_0(q)$ and $E_4(q)$ are connected by continuing analytically along the cycle ${\mathcal C}_q $ $(q=.224582 +.842777 I)$ (see Table 5). In section \ref{sec:pert} we obtained that the convergence radius of the series $E_4(q)$ is around $.875$. 
Hence for the eigenvalue $E_4(q)$ we also obtain compatibility.

For the anti-periodic case the closest branching point from the origin is $q=.535905+.640487 I$ ($|q|=.835115$) and the eigenvalues $E_1(q)$ and $E_3(q)$ are connected by continuing analytically along the cycle ${\mathcal C}_q $ $(q=.535905+.640487 I)$ (see Table 5). In section \ref{sec:pert} we obtained that the convergence radii of the series $E_1(q)$ and $E_3(q)$ are both around $.838$. For the eigenvalues $E_1(q)$ and $E_3 (q)$ we see compatibility and we obtain a reason why the convergence radii of $E_1(q)$ and $E_3(q)$ calculated in section \ref{sec:pert} are very close by considering the branching point.

We conclude that the convergence radii of the eigenvalues $E_m(q)$ $(m=0,1,2,3,4)$ calculated by perturbation and the locations of branching points calculated by considering the monodromy are compatible.

We presume that all eigenvalues $E_m(q)$ $(m\in 2\Zint _{\geq 0})$ in $\sigma _{{\bold H}_+} (H)$ (resp. all eigenvalues $E_m(q)$ $(m\in 2\Zint _{\geq 0}+1)$ in $\sigma _{{\bold H}_-} (H)$) are connected by analytic continuation in $q$.

{\bf Acknowledgments}

The author would like to thank Professor Hiroyuki Ochiai for valuable comments. 
Thanks are also due to the referee. He is partially supported by the Grant-in-Aid for Scientific Research (No. 15740108) from the Japan Society for the Promotion of Science.

\appendix
\section{Proof of Proposition \ref{prop:spec1}}
%\section{Proof of Proposition \ref{prop:spec1} and Proposition \ref{prop:spec2}}
To prove Proposition \ref{prop:spec1} we review some propositions from \cite{Tak1}, \cite{Tak3}.

Let ${\mathcal F}$ be the space spanned by meromorphic doubly periodic functions up to signs, namely
\begin{align}
& {\mathcal F}=\bigoplus _{\epsilon _1 , \epsilon _3 =\pm 1 } {\mathcal F} _{\epsilon _1 , \epsilon _3 }, \label{spaceF} \\
& {\mathcal F} _{\epsilon _1 , \epsilon _3 }=\{ f(x) \mbox{: meromorphic }| f(x+1)= \epsilon _1 f(x), \; f(x+\tau )= \epsilon _3 f(x) \} .
\end{align}
Let $V$ be the maximum finite-dimensional subspace in ${\mathcal F}$ which is invariant under the action of the Hamiltonian. Then it is known that $\dim V=2n+1$ \cite{WW}. Let $P(E)$ be the monic characteristic polynomial of the Hamiltonian $H$ (see Eq.(\ref{LameH})) on the space $V$, i.e. $P(E)=\prod_{i=1}^{n} (E-E_i)$ ($\{ E_i\} $ are eigenvalues of $H$ on $V$). Then the set $\sigma _{d} (H)$ coincides with the set of zeros of $P(E)$. From the periodicity we have $\sigma _d(H) \subset \sigma _{s} (H)$.

\begin{prop} \cite[Proposition 3.5]{Tak1}\label{prop:prod}
The equation
\begin{align}
& \left( \frac{d^3}{dx^3}-4\left( n(n+1)\wp (x)-E\right)\frac{d}{dx} -2\left(n(n+1)\wp '(x)\right) \right) \Xi (x,E)=0, \label{prodDE} 
\nonumber
\end{align}
has a nonzero doubly periodic solution which has the expansion
\begin{equation}
\Xi (x,E)=\sum_{j=0}^{n} b_j (E)\wp (x)^{n-j},
\label{Fx}
\end{equation}
where the coefficients $b_j(E)$ are polynomials in $E$, they do not have common divisors, and the polynomial $b_0(E)$ is monic.
Moreover the function $\Xi (x,E)$ is determined uniquely. 
\end{prop}

\begin{prop} \cite[Proposition 3.7]{Tak1}, \cite[Proposition 2.6]{Tak3} \label{prop:Lambda}
The function
\begin{equation}
\Lambda ( x, E)=\sqrt{\Xi (x,E)}\exp \int \frac{ \sqrt{-P(E)}dx}{\Xi (x,E)}
\label{integ1}
\end{equation}
is a solution to the differential equation (\ref{eqn:lame}).
\end{prop}
It follows from Eq.(\ref{integ1}) that, if $P(E) \neq 0$ then the functions $\Lambda (x, E)$ and $\Lambda (- x, E)$ are linearly independent (see also the proof of \cite[Lemma 3.6]{Tak1}) and they form a basis of the space of solutions to the differential equation (\ref{eqn:lame}).
Note that the function $\Lambda (x, E)$ is also expressed as
\begin{equation}
\Lambda (x, E) = A \prod _{i=1}^n \left( \frac{\sigma (x+a_i)}{\sigma (x)}e^{-x\zeta (a_i)} \right),
\end{equation}
for suitably chosen $A$ and $a_i$ $(i=1,\dots ,n)$ (see \cite[\S 39]{P} or \cite[\S 23.7]{WW}). 

From the periodicity of $\Xi (x,E)$ and the definition of $\Lambda (x, E)$, we have $\Lambda (x+1, E)=B(E) \Lambda (x,E)$ for some $B(E)$. 
Set $\Lambda^{\mbox{\scriptsize sym}}(x,E)=\Lambda (x,E)-(-1)^{n}\Lambda(-x,E)$. Then the relation $H\Lambda^{\mbox{\scriptsize sym}} (x,E) = E\Lambda ^{\mbox{\scriptsize sym}} (x,E)$ is obvious.

\begin{prop} (c.f. \cite[\S 4.4]{Tak3}) \label{prop:Tak344}
(i) If $B(E)=\pm 1$, then the function $\Lambda ^{\mbox{\scriptsize sym}} (x,E)$ is square-integrable on $(0,1)$. \\
(ii) If $P(E) \neq 0$, then the function  $\Lambda ^{\mbox{\scriptsize sym}} (x,E)$ is nonzero.\\
(iii) If $B(E)\neq \pm 1$, then $P(E)\neq 0$ and any nonzero solution to Eq.(\ref{eqn:lame}) is not square-integrable.
\end{prop}
\begin{proof}
(i) Because the exponents of the differential equation (\ref{eqn:lame}) at $x=0$ are $-n$ and $n+1$, we have the expansion $\Lambda ^{\mbox{\scriptsize sym}} (x,E)= x^{\alpha }(c_0+c_1x+ \cdots )$, where $c_0 \neq 0$ and ($\alpha =-n$ or $n+1$). From the property $\Lambda ^{\mbox{\scriptsize sym}} (x,E) =(-1)^{n+1} \Lambda ^{\mbox{\scriptsize sym}} (-x,E)$ and $n \in \Zint_{\geq 0}$, we have $\alpha =n+1$. Thus the function $\Lambda ^{\mbox{\scriptsize sym}} (x,E)$ is holomorphic at $x=0$.
It follows from the assumption $B(E)=\pm 1$ that $\Lambda ^{\mbox{\scriptsize sym}} (x+1,E) = \pm \Lambda ^{\mbox{\scriptsize sym}} (x,E)$. Hence the function  $\Lambda ^{\mbox{\scriptsize sym}} (x,E)$ is also holomorphic at $x=1$.
Since $\Lambda ^{\mbox{\scriptsize sym}} (x,E)$ satisfies the differential equation (\ref{eqn:lame}), it does not have singularity on the open interval $(0,1)$.
Therefore $\Lambda ^{\mbox{\scriptsize sym}} (x,E)$ is square-integrable on $(0,1)$.

(ii) It follows immediately from the linear independence of the functions $\Lambda (x,E)$ and $\Lambda(-x,E)$. 

(iii) Assume that $P(E)=0$. It follows from Eq.(\ref{integ1}) that $\Lambda (x,E)^2= \Xi (x,E)$. From the double-perioficity of the function $\Xi (x,E)$, we have $\Lambda (x+1,E)^2=\Lambda (x,E)^2$. Hence $\Lambda (x+1,E)=\pm \Lambda (x,E)$ and $B(E)=\pm 1$. Therefore we have $P(E)\neq 0$ under the assumption $B(E) \neq \pm 1$.

Assume that $B(E)\neq \pm 1$. Then $P(E)\neq 0$, and any solution to Eq.(\ref{eqn:lame}) can be written as a linear combination of $\Lambda (x,E)$ and $\Lambda (-x,E)$. The function $\Lambda (x, E)$ has poles at $x=0$ and $x=1$. Let $f(x)=C_1 \Lambda (x, E) +C_2 \Lambda (- x, E)$ be a non-zero square-integrable eigenfunction. The function $f(x)$ cannot have a pole at $x=0$ nor $x=1$ for square-integrability on $(0,1)$. If the function $f(x)$ is holomorphic at $x=0$, then we have $C_1 =-(-1)^n C_2$. From the periodicity we have $f(x+1)=C_1 B(E) \Lambda (x, E) +C_2 B(E)^{-1} \Lambda (- x, E) $. Hence we have $C_1 B(E)=-(-1)^n C_2 B(E)^{-1}$ for holomorphy of the function $f(x)$ at $x=1$. Under the assumption $B(E)\neq \pm 1$, we have $C_1=C_2=0$ and it contradicts to existence of the non-zero square-integrable eigenfunction.
\end{proof}

From Proposition \ref{prop:Tak344} (iii) we have $\sigma _{int} (H) \subset \sigma _{s} (H)$.
From Proposition \ref{prop:Tak344} (i), (ii)  we have $\sigma _{s} (H) \setminus \sigma _d(H) \subset \sigma _{int} (H)$.
Combining with $\sigma _{d} (H) \subset \sigma _{s} (H)$ we have $\sigma _{s} (H) = \sigma _d(H) \cup  \sigma _{int} (H)$.
Therefore we obtain Proposition \ref{prop:spec1} (i). To prove Proposition \ref{prop:spec1} (ii), it is sufficient to show the following lemma:
\begin{lemma} \label{lemma:intersect}
If $q (=\exp(\pi \sqrt{-1} \tau)) \in \Rea$ and $0<|q|<1$, then
%If $-1<q(=\exp (\pi\sqrt{-1}\tau))<1$, then
\begin{equation}
\sigma _d(H) \cap  \sigma _{int} (H) = \phi .
\end{equation}
\end{lemma}
\begin{proof}
By Proposition \ref{prop:spec2}, it is enough to show that $\sigma _d(H) \cap  \sigma _{{\bold H}} (H) = \phi$. Set
\begin{align}
& I = -\frac{d ^2}{d x^2} +n(n+1) \wp (x+\tau /2). 
\end{align}
Then the potential does not have poles on $\Rea $. 
Set
\begin{equation}
\tilde{{\bold H}} = \left\{ f \! : \Rea \rightarrow \Cplx \mbox{ measurable} \left|
\begin{array}{l}
 \int_{0}^{1} |f(x)| ^2 dx<+\infty, \\
 f(x)= f(x+2) \mbox{ a.e. }x, \\
\end{array}
\right. \right\}
\end{equation}
From the periodicity we have $\sigma _s (H) = \sigma  _{\tilde{\bold H}} (I)$ as a set. 

If $q=0$ then a basis of eigenfunctions in the space $\tilde{\bold H}$ is $\{ \exp( m \pi \sqrt{-1} x ) \} _{m \in \Zint }$ and we have $\sigma  _{\tilde{\bold H}} (I) = \{ \pi ^2 m^2 -\pi ^2 n(n+1)/3 | \; m \in \Zint \}$ with multiplicity.  For the case $q=0$, the set $\{ \pi ^2 m^2 -\pi ^2 n(n+1)/3 | \; m \in \Zint, \;m \geq  n+1 \}$ coincides with the set $\sigma _s (H)$.
The set $\sigma _d (H)$ tends to the set $\{ \pi ^2 m^2 -\pi ^2 n(n+1)/3 | \; m \in \Zint, \; -n\leq m \leq n  \}$ as $q \rightarrow 0$.
We define the set $\sigma _d (H)$ for the case $q=0$ by $\sigma _d (H) = \{ \pi ^2 m^2 -\pi ^2 n(n+1)/3 | \; m \in \Zint, \; -n\leq m \leq n  \}$.
Then we can check directly that $\sigma _s (H) = \sigma _d(H) \cup  \sigma _{int} (H) $ for the case $q=0$.

By a similar discussion to Proposition \ref{mainthmKato} and \cite[Proposition 3.3]{Tak2} (see also \cite{Kat}), it follows that all eigenvalues of $I$ $(-1<q<1)$ on the space $\tilde{\bf H}$ can be represented as $\tilde{E}_m (q)$ $(m \in \Zint )$, which is real-holomorphic in $q \in (-1,1)$, $\tilde{E}_m (0) = \pi ^2 m^2 -\pi ^2 n(n+1)/3 $ and the operator $I$ $(-1<q<1)$ forms a holomorphic family of type (A) (for definition see \cite{Kat}).
From the equation $\sigma _s (H) = \sigma _d(H) \cup  \sigma _{int} (H) = \sigma _d(H) \cup  \sigma _{{\bold H}} (H) =\sigma  _{\tilde{\bold H}} (I)$ and that elements in $\sigma _d (H)$, $\sigma _{{\bold H}}(H)$ and $\sigma _{\tilde{\bold H}} (I)$ are all real-holomorphic in $q$ $(-1<q<1)$, we have $\sigma _d (H)= \{ \tilde{E}_m (q) | \;  m \in \Zint, \; -n\leq m \leq n  \}$ and $\sigma _{{\bold H}} (H) =\{ \tilde{E}_m (q) | \;  m \in \Zint, \; m \geq n+1  \}$. Moreover we have $\tilde{E}_{m+n+1} (q) = E_m (q)$ $(m \in \Zint _{\geq 0})$ and the multiplicity of the eigenvalue $\tilde{E}_{m+n+1} (q)$ $(m \in \Zint _{\geq 0})$ on the space $\tilde{\bold H}$ is two.

Suppose $E \in \sigma _d(H) \cap  \sigma _{{\bold H}} (H)$. Then $E$ is both the eigenvalue in $\sigma _d(H) \subset \sigma _{\tilde{\bold H}} (I)$ (multiplicity $\geq 1$) and the eigenvalue in $\sigma _{{\bold H}} (H) \subset \sigma _{\tilde{\bold H}} (I)$ (multiplicity $\geq 2$) and the multiplicity is summed up because the operator $I$ $(-1<q<1)$ form a holomorphic family of type (A). Hence the multiplicity of the eigenvalue $E$ is no less than three. However that is impossible because the dimension of the solution to the second-order linear ordinary differential equation $(I-E)f(x)=0$ with the boundary condition $f(x) \in \tilde{\bf H}$ is no more than two.
Thus we obtain that if $-1<q<1 $ then $\sigma _d(H) \cap  \sigma _{{\bold H}} (H)= \phi$.
\end{proof}

\section{}
We note definitions and formulas for elliptic functions.
Let $\omega_1$ and $\omega_3$ be complex numbers such that the value $\omega_3/ \omega_1$ is an element of the upper half plane.

The Weierstrass $\wp$-function, the Weierstrass sigma-function and the Weierstrass zeta-function are defined as follows:
\begin{align}
& \wp (x)=\wp(x|2\omega_1, 2\omega_3)= \\
& \; \; \; \; \frac{1}{x^2}+\sum_{(m,n)\in \Zint \times \Zint \setminus \{ (0,0)\} } \left( \frac{1}{(x-2m\omega_1 -2n\omega_3)^2}-\frac{1}{(2m\omega_1 +2n\omega_3)^2}\right), \nonumber \\
& \sigma (x)=\sigma (x|2\omega_1, 2\omega_3)=x\prod_{(m,n)\in \Zint \times \Zint \setminus \{(0,0)\} } \left(1-\frac{x}{2m\omega_1 +2n\omega_3}\right) \cdot \nonumber \\
& \; \; \; \;\cdot \exp\left(\frac{x}{2m\omega_1 +2n\omega_3}+\frac{x^2}{2(2m\omega_1 +2n\omega_3)^2}\right), \nonumber \\
& \zeta(x)=\frac{\sigma'(x)}{\sigma (x)}. \nonumber
\end{align}
Setting $\omega_2=-\omega_1-\omega_3$ and 
\begin{align}
& e_i=\wp(\omega_i), \; \; \; \eta_i=\zeta(\omega_i), \; \; \; \; (i=1,2,3) .
\end{align}
yields the relations
\begin{align}
& e_1+e_2+e_3=\eta_1+\eta_2+\eta_3=0, \; \; \; \wp(x)=-\zeta'(x), \label{eq:sigmaper} \\
& (\wp'(x))^2=4(\wp(x)-e_1)(\wp(x)-e_2)(\wp(x)-e_3), \nonumber \\
& \zeta(x+2\omega_i)=\zeta(x)+2\eta_i, \; \; \; \sigma(x+2\omega_i)=-\sigma(x)e^{2\eta_i(x+\omega _i)},  \nonumber \\
& \wp(x+2\omega_i)=\wp(x),  \; \; \; \wp (x+\omega _i)-e_i= \frac{(e_i-e_{i'})(e_i-e_{i''})}{\wp (x)-e_i} , \nonumber
\end{align}
where $\{i, i', i'' \} =\{1,2,3\}$.
On elliptic integrals we have
\begin{align}
& t-\omega _i= \int _{e_i}^{\wp (t)} \frac{ds}{\sqrt{4(s-e_1)(s-e_2)(s-e_3)}}, \label{ellinteg1} \\
& \zeta (t)-\eta _i= \int _{e_i}^{\wp (t)} \frac{-sds}{\sqrt{4(s-e_1)(s-e_2)(s-e_3)}}, \; \; \; (i=1,2,3). \nonumber
\end{align}
The co-$\wp $ functions $\wp_i(x)$ $(i=1,2,3)$ are defined by
\begin{align}
& \wp_i(x)=\exp (-\eta_i x)\sigma(x+\omega_i)/( \sigma(x)\sigma(\omega _i)),
\label{wpi}
\end{align}
and satisfy
\begin{align}
& \wp _i(x) ^2 =\wp(x)-e_i,\; \; \; \; \; \; \; (i=1,2,3). \label{sigmai}
\end{align}

Set $\omega_1=1/2$, $\omega_3=\tau /2$ and $q=\exp (\pi \sqrt{-1} \tau )$. The expansions of the Weierstrass $\wp$ function, the Weierstrass $\zeta$ function and $\eta_1$ in the variable $q$ are written as follows: 
\begin{align}
& \wp (x)=
-2\eta_1 +\frac{\pi^2 }{\sin^2 (\pi x)} -8\pi^2 \sum_{k=1}^{\infty} \frac{kq^{2k}}{1-q ^{2k}} \cos 2k \pi x , \label{wpth} \\
& \zeta (x)=
2\eta_1x +\frac{\pi }{\tan (\pi x)} +4\pi \sum_{k=1}^{\infty} \frac{q ^{2k}}{1-q ^{2k}} \sin 2k \pi x , \nonumber \\
& \eta_1 = \pi^2 \left( \frac{1}{6} - 4\sum_{k=1}^{\infty} \frac{kq ^{2k}}{1-q ^{2k}} \right) . \nonumber 
\end{align}

\end{document}